\documentclass[11pt,  reqno]
{amsart}

\usepackage{amsmath,amssymb,amscd,amsthm,amsxtra, esint}

\allowdisplaybreaks[2]

\newtheorem{theorem}{Theorem} [section]

\newtheorem{lemma}[theorem]{Lemma}

\newtheorem{remark}[theorem]{Remark}

\newtheorem{definition}[theorem]{Definition}

\newtheorem*{acknowledgment}{Acknowledgment}

\newcommand{\noi}{\noindent}
\newcommand{\Z}{\mathbb{Z}}
\newcommand{\R}{\mathbb{R}}

\newcommand{\T}{\mathbb{T}}

\renewcommand{\L}{\mathcal{L}}

\newcommand{\N}{\mathcal{N}}

\newcommand{\A}{\mathcal{A}}

\newcommand{\al}{\alpha}

\newcommand{\dl}{\delta}

\newcommand{\nb}{\nabla}

\newcommand{\eps}{\varepsilon}

\newcommand{\g}{\gamma}
\newcommand{\G}{\Gamma}
\newcommand{\ld}{\lambda}

\newcommand{\s}{\sigma}

\newcommand{\ft}{\widehat}

\newcommand{\wt}{\widetilde}
\newcommand{\cj}{\overline}
\newcommand{\dx}{\partial_x}
\newcommand{\dt}{\partial_t}

\renewcommand{\l}{\ell}
\renewcommand{\o}{\omega}

\newcommand{\les}{\lesssim}
\newcommand{\ges}{\gtrsim}

\newcommand{\jb}[1]
{\langle #1 \rangle}

\newcommand{\too}{\longrightarrow}

\newcommand{\oo}{\pmb{\o}}
\newcommand{\n}{\pmb{n}}

\numberwithin{equation}{section}
\numberwithin{theorem}{section}

\begin{document}

\baselineskip = 14pt

\title[Defocusing NLS with limit periodic initial data]
{Global existence for  the 
defocusing nonlinear Schr\"odinger equations
with limit periodic initial data}

\dedicatory{
Dedicated to Professor Gustavo Ponce on the
occasion of his sixtieth birthday}

\author{Tadahiro Oh}

\address{
Tadahiro Oh\\
School of Mathematics\\
The University of Edinburgh, 
and The Maxwell Institute for the Mathematical Sciences\\
James Clerk Maxwell Building\\
The King's Buildings\\ 
Peter Guthrie Tait Road\\
Edinburgh\\ 
EH9 3FD\\ United Kingdom}

\email{hiro.oh@ed.ac.uk}

%

\begin{abstract}

We consider the Cauchy problem 
for the defocusing nonlinear Schr\"odinger equations (NLS)  on the real line
with a special subclass of 
almost periodic functions as initial data.
In particular, 
we prove global existence of solutions to NLS
with limit periodic functions as initial data
under some regularity assumption.

\end{abstract}



\maketitle

\section{Introduction}
\label{SEC:1}

We consider the Cauchy problem 
for the  defocusing  nonlinear Schr\"odinger equations (NLS) on $\R$:
\begin{align}
\begin{cases}
i \dt u + \dx^2 u =  |u|^{2k} u, \\
u|_{t = 0} = f,
\end{cases}
\quad (t, x) \in \R\times \R, 
\label{NLS1}
\end{align}

\noi
for $k \in \mathbb{N}$.
The Cauchy problem \eqref{NLS1} has been studied extensively
in terms of the usual 
 Sobolev spaces $H^s(\R)$ on the real line and 
 the Sobolev spaces $H^s_{\text{per}}(\R) \simeq H^s(\T)$,
 $\T = \R/\Z$,  
 of periodic functions (of a fixed period) on $\R$.
See, for example,  Ginibre-Velo \cite{GV}, Tsutsumi \cite{Tsutsumi},
and Bourgain \cite{BO}.
There are several known conservation laws for NLS \eqref{NLS1}. 
In particular, %
the conservation 
of the Hamiltonian and the mass, defined by 
\begin{align*}
\text{Hamiltonian: } & H[u](t) = \frac{1}{2}\int_{\mathcal{M}} |\dx u(t, x)|^2 dx
+ \frac{1}{2k+2} \int_{\mathcal{M}} |u(t, x)|^{2k+2} dx,  \\
\text{Mass: } & Q[u](t) = \int_{\mathcal{M}} |u(t, x)|^2 dx, 
\end{align*}

\noi
where $\mathcal{M} = \R$ or $\T$, 
plays a crucial role
in establishing global well-posedness of \eqref{NLS1}.
See also  the monographs \cite{CAZ, TAO} for more references on the subject.

Our main interest in this paper is to study 
global-in-time  behavior of solutions to the Cauchy problem \eqref{NLS1}
with a particular subclass of  almost periodic functions as initial data.
In particular, 
we prove global existence of unique solutions 
to  \eqref{NLS1} with {\it limit periodic} functions
(see Definitions \ref{DEF:LP} and \ref{DEF:LP2} below)
as initial data
under some regularity assumption.

In Subsection \ref{SUBSEC:AP}, we go over the basic definitions
and properties of almost periodic functions
along with 
the known well-posedness results
for NLS \eqref{NLS1}
with almost periodic functions as initial data.
We then introduce limit periodic functions
and state
our main result (Theorem \ref{THM:1})
in Subsection \ref{SUBSEC:LP}.

\subsection{Almost periodic functions}
\label{SUBSEC:AP}

Let us first recall the definition of 
almost periodic functions  due to  Bohr \cite{Bohr}.

\begin{definition}\label{DEF:1}\rm 
We say that a 
function $f$ on $\R$ is  {\it almost periodic},
if it is continuous and,  for every $\eps > 0$, there exists  $L = L(\eps, f)>0$
such that every interval of length $L$ on $\R$ contains
 a number $\tau$ such that 
\begin{align*}
 \sup_{x \in \R} |f(x-\tau) - f(x)|< \eps.
 \end{align*}

\noi
We use $AP(\R)$ to denote the space of almost periodic functions on $\R$.

\end{definition}

\noi
In the following, we briefly go over some basic properties of almost periodic functions.
See Besicovitch \cite{B}, Corduneanu \cite{C}, and Katznelson \cite{K}
for 
more on the subject.
It is well known that the following two notions in Definition \ref{DEF:2} (i) and (ii) are
equivalent to the notion of 
 almost periodic functions in Definition \ref{DEF:1}.

\begin{definition}\label{DEF:2}\rm

\noi
(i) We say that a function $f$ on $\R$ has the {\it approximation property},
if it can be uniformly approximated by trigonometric polynomials
(of finite degrees).

\noi
(ii) We say that a continuous function on $\R$ is {\it normal}
if, given any sequence $\{x_n\}_{n = 1}^\infty \subset \R$,
the collection $\{ f(\,\cdot  + x_n)\}_{n = 1}^\infty$ is precompact in $L^\infty(\R)$.

\end{definition}

%

\begin{remark}\rm
Definition \ref{DEF:2} (i) states that 
the collection $AP(\R)$ of almost periodic functions is exactly the closure
of the trigonometric polynomials
with respect to the uniform metric induced by the $L^\infty$-norm.
There are also the notions of different classes of generalized almost periodic functions due to Stepanov,
Weyl, and Besicovitch
by considering the closures of the trigonometric polynomials
under different metrics.
 The corresponding spaces are denoted by $S^p$,  $W^p$,  and $B^p$, respectively.
Note that  we have $AP(\R) \subset S^p \subset W^p \subset B^p$, $p \geq 1$.
For example, see  
 Remarks \ref{REM:BP} and \ref{REM:SP} below.
In the literature, 
almost periodic functions in $AP(\R)$
are sometimes
referred to as {\it uniformly} almost periodic functions
in order to distinguish them from these generalized almost periodic functions.
In this paper, however, we only consider almost periodic functions
in Bohr's sense
according to Definition \ref{DEF:1}.

\end{remark}

The space $AP(\R)$ of almost periodic functions 
is a closed subalgebra of $L^\infty(\R)$
and 
almost periodic functions are uniformly continuous.
 Given $f \in AP(\R)$,
 we can define the so-called
mean value $M(f)$  of $f$
 by
\begin{align*}
	 M(f) := \lim_{L\to \infty} \frac{1}{2L}\int_{-L}^L f(x) dx.
\end{align*}

\noi
We then define an inner product $\jb{\cdot, \cdot}_{\L^2}$ on $AP(\R)$ by
\begin{align*}
\jb{f, g}_{\L^2} := M(f\cj{g}) = \lim_{L\to \infty} \frac{1}{2L}\int_{-L}^L f(x) \cj{g}(x) dx
\end{align*}

\noi	
for $f, g \in AP(\R)$.
Note that $M(f)$ and $\jb{f, g}_{\L^2}$
are  well defined for $f, g \in AP(\R)$.
This inner product induces
 the $\L^2$-norm defined by 
\begin{align}
\|f\|_{\L^2} := M(|f|^2) 
= \jb{f, f}_{\L^2}
\label{AP2}
\end{align}

\noi
and makes the space $AP(\R)$ of almost periodic functions
 a pre-Hilbert space (missing completeness).
Note that the $\L^2$-norm is a norm on $AP(\R)$, 
but not for general functions on $\R$.
For example, we have $\|f\|_{\L^2} = 0$ for any bounded function $f$ with a compact support.

\begin{remark} \label{REM:BP}\rm
The space $B^p$, $p \geq 1$,  of Besicovitch's generalized almost periodic functions
is precisely the closure of the trigonometric polynomials
under the $B^p$-metric $d_{B^p}$ defined by 
\[ d_{B^p}(f, g) 
:= \limsup_{L\to \infty} \bigg(\frac{1}{2L}\int_{-L}^L |f(x)- g(x) |^p dx\bigg)^\frac{1}{p}.\]

\noi
It is known that $B^2$ is complete with respect to the
$\mathcal{L}^2$-norm defined in \eqref{AP2}.
See \cite[Riesc-Fischer Theorem on p.\,109]{B}.

\end{remark}

Now, let us turn our attention to the Fourier analysis of almost periodic functions.
The complex exponentials
$\{ e^{ 2 \pi i \o x}\}_{\o \in \R}$ form an orthonormal family
under the inner product $\jb{\cdot, \cdot}_{\mathcal{L}^2}$.
Given $f \in AP(\R)$, 
we then define its Fourier coefficient  by
\begin{align}
\ft f (\o) 
: = \jb {f, e^{2\pi  i \o x}}_{\L^2} = M(f e^{- 2 \pi i \o x}).
\label{AP2a}
\end{align}

\noi
It follows from Bessel's inequality that 
$\ft f(\o) = 0$ except
for  countable many values of $\o \in \R$.
Given $f \in AP(\R)$, 
we define the Fourier series associated to $f$ by 
\begin{equation}
 f (x) \sim \sum_{\o \in \s(f)} \ft f (\o) e^{2\pi i \o x}.
\label{AP3c}
 \end{equation}

\noi
Here, 
$\s(f)$ denotes the frequency set of $f$ defined by 
\[\s(f) : = \{ \o \in \R:\, \ft f(\o) \ne 0\}.\]

\noi
It is known that the orthonormal family
$\{ e^{ 2\pi i \o x}\}_{\o \in \R}$ is complete
in the sense that two distinct almost periodic functions
have distinct Fourier series.
Moreover, we have the Parseval's identity:
\begin{align*}
  \|f\|_{\L^2} = \bigg(\sum_{\o \in \R} |\ft f (\o)|^2 \bigg)^\frac{1}{2}
  \end{align*}

\noi
for $f \in AP(\R)$.
Regarding the actual convergence of the Fourier series
to an almost periodic function, we have the following lemma.

\begin{lemma}[Theorem 1.20 in \cite{C}] \label{LEM:AP2}
Let $f \in AP(\R)$.
If the Fourier series associated to $f$ converges uniformly on $\R$,
then it converges to $f$.  Namely, we have
\begin{align}
 f (x) = \sum_{\o \in \s(f)} \ft f (\o) e^{2\pi i \o x}.
\label{AP3a}
 \end{align}
\end{lemma}

\noi
Since our argument in this paper is based on the Fourier series representation
\eqref{AP3a}, 
Lemma \ref{LEM:AP2} plays an important role in the following.

Given $\oo = \{\o_j\}_{j = 1}^\infty
\in \R^{\mathbb N}$,  we
say that $\oo$ is linear independent
if any relation of the form:
\[ \sum_{j = 1}^N \al_j \o_j =0, \quad \al_j  \in \mathbb Q, \]

\noi
implies that $\al_j = 0$, $j = 1, \dots, N$.

Given a set $S$ of real numbers,
we say that a linearly independent
set $\oo = \{\o_j\}_{j = 1}^\infty$
is a basis for the set $S$, if every element in $ S$
can be represented as
a finite linear combination of elements in $\oo$
with rational coefficients.
Given $f \in AP(\R)$,
we say that a linearly independent set  $\pmb{\o} = \{\o_j\}_{j = 1}^N$,
allowing the case $N = \infty$,
is a basis of $f$, if
it is a  basis for
the frequency set
$\s(f)$ of $f$.
Lemma 1.14 in \cite{C} guarantees existence of a basis of $f \in AP(\R)$.
We say that a basis $\pmb{\o} = \{\o_j\}_{j = 1}^N$ of $f$ is an integral basis
if any element in the frequency set $\s(f)$ can be written as a finite
linear combination of elements in $\oo$
with integer coefficients.
If there exists a finite integral basis of $f$, i.e.~$N < \infty$, then we say that the function $f$
is {\it quasi-periodic}.

We conclude this subsection 
by going over 
 the known well-posedness results
on NLS with almost periodic functions as initial data.
In the periodic setting, 
Bourgain \cite{BO} proved
local well-posedness of \eqref{NLS1}
in 
$H^s(\T)$
with $s = 0$ when $k = 1$
and 
with 
$s > \frac{1}{2} - \frac 1k$ when $k \geq 2$.
The conservation of the Hamiltonian and the mass
then yields global well-posedness of \eqref{NLS1}
in $L^2(\T)$ when $ k = 1$
and in $H^1(\T)$ when $ k\geq 2$.

Regarding local well-posedness with quasi-periodic initial data, 
Tsugawa  \cite{Tsugawa} 
considered quasi-periodic functions $f$ of the form:
\begin{align}
	 f(x) = \sum_{\n \in \Z^N} \ft f(\pmb{\o}\cdot \pmb{n}) e^{2 \pi i (\oo \cdot \n)x},
\label{AP4a}
\end{align}

\noi
for  a frequency set $\oo = \{\o_j\}_{j = 1}^N \in \R^N$ with $N \in \mathbb N$.
Then, he proved
local well-posedness of the Korteweg-de Vries equation (KdV) on $\R$:
\begin{align*}
\dt u + \dx^3 u = u  \dx u
\end{align*}

\noi
with quasi-periodic initial data of the form \eqref{AP4a}
under some regularity condition.
Moreover, 
by defining 
in a Sobolev-type space $\mathcal{H}^{\pmb{s}}_{\oo}(\R)$
of quasi-periodic functions of the form \eqref{AP4a} by the norm
 \begin{align*}
  \|f\|_{\mathcal{H}^{\pmb{s}}_{\oo}(\R)}
 :=  \|    \jb{\n}^{\pmb{s}}\ft f (\oo\cdot \n) \|_{\l^2_{\n}(\Z^N)},
 \quad
 \jb{\n}^{\pmb{s}} : =  \prod_{j = 1}^N \big( 1+ |n_j|^2\big)^\frac{s_j}{2},
 \end{align*}

\noi
where  $\pmb{s} = \{s_j\}_{j = 1}^N \in \R^{N}$, 
\cite[Lemma 2.2 (i)]{Tsugawa}
implies that NLS \eqref{NLS1} is locally well-posed
in  $\mathcal{H}^{\pmb{s}}_{\oo}(\R)$, 
provided that $\min(s_1, \dots, s_N) > \frac{1}{2}$.

Let us turn our attention
to the  generic almost periodic setting. 
Fix  a frequency set $\pmb{\o} = \{\o_j\}_{j = 1}^\infty \in \R^{\mathbb N}$.
For almost periodic functions $f \in AP(\R)$
of the form:
\begin{align}
	 f(x) \sim \sum_{\n \in \Z^\mathbb{N}} \ft f(\pmb{\o}\cdot \pmb{n}) e^{2\pi i (\oo \cdot \n)x},
\label{AP4}
\end{align}

\noi
where $\n = \{ n_j\}_{j = 1}^\infty \in \Z^\mathbb{N}$, 
we  define 
 the $\A_{\oo}$-norm  by
\begin{align*}
\| f\|_{\A_{\oo}(\R)} = \| \ft f(\oo\cdot \n)\|_{\l^1_{\n}(\Z^\mathbb{N})}.
\end{align*}

\noi
Then, we define 
the algebra $\A_{\oo}(\R)$ by
\begin{align}
\A_{\oo}(\R) & = \big\{ f \in AP(\R):\,  \notag \\
& \hphantom{XXXXX}
f \text{ is of the form } \eqref{AP4} \text{ and }
 \| f\|_{\A_{\oo}(\R)} < \infty\big\}.
\label{AP44}
 \end{align}

\noi

\noi
See \cite{Oh} for some basic properties of $\A_{\oo}(\R)$.
In \cite{Oh}, we proved local well-posedness
in $\A_{\oo}(\R)$
for  NLS with a power-type nonlinearity, including \eqref{NLS1}.

In view of the result in \cite{Oh}, 
it is natural to consider the global-in-time behavior
of solutions to \eqref{NLS1}.
This is an extremely difficult question in general.
There are, however, several known
global existence results for  the cubic NLS,  \eqref{NLS1} with $k = 1$,
and KdV in the almost periodic and quasi-periodic settings.
Egorova \cite{E} and Boutet de Monvel-Egorova \cite{BE}
constructed global-in-time solutions
to KdV and  the cubic NLS with almost periodic initial data,
assuming some conditions,
including Cantor-like spectra for
the  corresponding Schr\"odinger operator (for KdV)
and Dirac operator (for the cubic NLS).
In particular,
the class of almost periodic initial data in \cite{E, BE}
includes
almost periodic functions $f$
that can be approximated by periodic functions $f_j$
of growing periods $\al_j \to \infty$
in $S^{s, 2}(\R)$
with $s \geq 4$ for KdV and $s \geq 3$ for cubic NLS.
See \eqref{Stepanov2} below
for the definition of the $S^{s, 2}$-norm.
Moreover, the  convergence of $f_j$ to $f$ in 
the $S^{s, 2}$-norm  is assumed to be exponentially fast.
It is worthwhile to mention that the solutions constructed in \cite{E, BE}
are almost periodic in both $t$ and $x$.
There is also a recent global existence result of KdV
with quasi-periodic initial data
by Damanik-Goldstein \cite{DG}.
Their result states that if  the Fourier coefficient $\ft f(\oo\cdot \n)$
of a quasi-periodic initial condition \eqref{AP4a}
decays exponentially fast (in $\n$),
then there exists a unique global solution
whose Fourier coefficient also decays exponentially fast
(with a slightly worse constant),
provided that a smallness condition on the initial condition $f$
and a Diophantine condition on $\oo$ are satisfied.

We emphasize that
the above results rely heavily on the inverse spectral method and
on the complete integrability of the equations.
These methods are not applicable to the non-integrable case,
i.e.~ \eqref{NLS1} with $k \geq 2$.
Our main goal in this paper is to 
establish global existence 
via an analytical method
without complete integrability.

\subsection{Limit periodic functions
and the 
main result}
\label{SUBSEC:LP}

In the following, we restrict our attention
to a particular subclass of almost periodic functions, 
called limit periodic functions.

\begin{definition}\label{DEF:LP}\rm
We say that a function $f$ on $\R$ is {\it limit periodic}
if it is a uniform limit
of 
continuous periodic functions.
\end{definition}

\noi
Note that  a limit periodic function is almost periodic,
since 
$AP(\R)$ is closed under the $L^\infty$-norm.
The following characterization of limit period functions plays
an essential role in our analysis.

\begin{lemma}[Theorem and Converse Theorem on p.\,32 in \cite{B}]\label{LEM:PL}
 An almost periodic function $f$ on $\R$ is limit periodic
if and only if its Fourier series is given by 
\begin{align}
 f(x) \sim \sum_{m = 1}^\infty \ft f(r_m \o) e^{ 2\pi i r_m \o x}
\label{LP1}
 \end{align}

\noi
for some $\o \in \R$ and $\{ r_m\}_{m = 1}^\infty  \subset  \mathbb Q$.
\end{lemma}

\noi
Namely, an almost periodic function is
limit periodic if and only if it has a one-term basis, 
i.e.~all the frequencies are rational multiples of a single frequency $\o \in \R$.

In view of the global well-posedness result \cite{BO} in the periodic setting, 
we assume that our initial condition $f$ is not periodic.
Note that 
 an almost periodic function is periodic
if and only if it has a one-term  {\it integral}  basis. 
Hence, we assume that $\o$ 
(or any of its rational multiple) is not an integral basis.\footnote{This in particular implies
that the denominators of $\{r_m\}_{m \in \mathbb{\N}}$ are unbounded.} %
It also follows from Lemma \ref{LEM:PL}
that if a limit periodic function is quasi-periodic, 
then it is periodic.
Therefore, we consider limit periodic functions 
that are not quasi-periodic in the following.
Lastly, 
we also assume that $\o \ne 0$
 in the following,
since $\o= 0$  corresponds to constant functions.

\begin{definition} \label{DEF:LP2}
\rm 
Let $\o \in \R \setminus \{0\}$.
We denote the class of limit periodic functions with 
a one-term basis $\o\in \R$ by $LP_\o(\R)$.
Namely, 
we say that $f \in LP_\o(\R)$ 
if it has the Fourier series expansion \eqref{LP1}
with this specific $\o$.

\end{definition}

Given a limit periodic function $f$, define $f_j$,  $j \in \mathbb N$,  by 
\begin{align}
 f_j(x) \sim \sum_{ m \in A(j)} \ft f(r_m\o) e^{ 2\pi i (r_m j!) \, \frac{ \o}{j!} x}, 
\label{LP2}
 \end{align}

\noi
where $A(j)$ is given by 
\begin{equation}
A(j) = \big\{ m \in \mathbb N:\, r_m j! \in \Z \big\}.
\label{LP2a}
\end{equation}

\noi
Then, $f_j$ is periodic with period 
\begin{equation}
L_j : = \frac{j! }{ \o} .
\label{period}
\end{equation}

\noi
Moreover, it is known that 
$f_j$ converges to $f$ uniformly as $j \to \infty$.
See \cite[p.45]{B}.
We refer to $f_j$ as the periodization of $f$
(with period $L_j$).
See Lemma \ref{LEM:ave} below.

Let $p \geq 1$ and $s \in \mathbb{N}$. 
We define our function space $S^{s, p}(\R)$ 
by 
\[ S^{s, p}(\R) = \{ f \in L^1_\text{loc}(\R):\,  \| f\|_{S^{s, p}(\R)} < \infty\},\]

\noi
where
the $S^{s, p}$-norm is defined by 
\begin{align}
\| f\|_{S^{s, p}(\R)}
:= \sup_{y\in \R} \bigg(\int_{y}^{y+1} |f(x) |^p + |\dx^s f(x) |^pdx \bigg)^\frac{1}{p}.
\label{Stepanov2}
\end{align}

\noi
Note that $S^{s, p}(\R)$ is complete,\footnote{To see this, 
we can consider the following $\wt S^{s, p}$-norm given by 
\begin{align*}
\| f\|_{\wt S^{s, p}(\R)}
:= \sup_{j\in \Z} \bigg(\int_{j}^{j+1} |f(x) |^p + |\dx^s f(x) |^pdx \bigg)^\frac{1}{p}.
\end{align*}

\noi
Clearly, the $S^{s, p}$- and $\wt S^{s, p}$-norms are equivalent.
Moreover,  $\wt S^{s, p}(\R)$ is complete, since 
$W^{s, p}([j, j+1))$ is complete for each $j \in \Z$.
} just like the usual Sobolev spaces $W^{s, p}([y, y+1))$.
On the one hand,  the definition of $S^{s, p}(\R)$ has nothing to do with almost periodic functions.
On the other hand, we point out that 
the $S^{s, p}$-norm is closely related to the $S^p$-metric used for
Stepanov's generalized almost periodic functions.
See Remark \ref{REM:SP} below.

Now, we are ready to state our main result.

\begin{theorem}\label{THM:1}
Given $\o \in \R  \setminus \{0\} $, let $f \in LP_\o(\R) \cap S^{1, 2}(\R)$
and $f_j$ be the periodization of $f$ as in \eqref{LP2}.

\smallskip
\noi
\textup{(i)}
Suppose that 
  there exist  $\eps > 0$ and $B > 0$ such that 
\begin{align}
e^{L_{j+1}^{2k + \eps}} \| f _j - f\|_{S^{1, 2}(\R)} \leq B 
\label{LP3}
\end{align}

\noi	
for all sufficiently large $j$.
Then, there exists a unique global solution $u 
\in C(\R; S^{1, 2}(\R)) \subset C(\R; L^\infty(\R))$
 to the defocusing NLS \eqref{NLS1}
with $u|_{t = 0 } = f$.
Moreover, $u(t)$ lies in $LP_\o(\R)$ for each $t \in \R$.

\smallskip
\noi
\textup{(ii)}
Given $J\in \mathbb{N}$ and $K > 0$, 
define $\mathcal{B}^\o(J, K)$ by 
\begin{align}
\mathcal{B}^\o (J, K)
& = \big\{ f \in LP_\o(\R) \cap S^{1, 2}(\R): \, 
\|f\|_{S^{1, 2}(\R)} \leq K, \notag \\
& \hphantom{XXXXXXXX}
\eqref{LP3}
\textup{ is satisfied for all }j \geq J\big\}.
\label{BJ}
\end{align}

\noi
Then, for fixed $t \in \R$, 
the solution map$\, : u(0) = f \mapsto u(t)$ constructed in \textup{(i)}
is continuous on $\mathcal{B}^\o(J, K)$
with the $S^{1, 2}$-topology.
\end{theorem}

\noi
To the best of the author's knowledge, 
Theorem \ref{THM:1}
is the first global existence result
for NLS \eqref{NLS1}, $k \geq 2$, 
with limit periodic functions as initial data.
While global existence of solutions for the cubic NLS ($k = 1$)
was previously proved in  \cite{BE}, 
the argument in \cite{BE} relied heavily on the complete integrability 
of the cubic NLS
and is not applicable to the non-integrable case $k \geq 2$.
We prove  Theorem \ref{THM:1}
by combining global well-posedness
of the defocusing  NLS in the periodic setting
and scaling invariance.

\begin{remark}\rm 
(i)
Given $f \in LP_\o(\R)$, 
let $r_m$ be as in \eqref{LP1}.
Then, under the hypothesis of Theorem \ref{THM:1}, 
we prove that
$u(t) \in A_{\oo}(\R)$ for each $t\in \R$,
where 
\[\oo := \{ r_m \o \}_{m = 1}^\infty \in \R^\mathbb{N}.\]

\noi 
Then, from 
the local well-posedness  in 
$\A_{\oo}(\R)$ \cite{Oh}, 
we obtain uniqueness and  (local-in-time) continuous dependence,  at each $t \in \R$, 
in the  $\A_{\oo}(\R)$-topology
of the global-in-time flow  constructed in Theorem \ref{THM:1}.

\smallskip

\noi
(ii)  The uniqueness statement in Theorem \ref{THM:1} holds in $C(\R; \A_{\oo}(\R))$
as mentioned above.
There is also a mild uniqueness statement\footnote{This is analogous to having uniqueness only as a limit of classical solutions to some evolution equations (even in the usual Sobolev spaces).}  as a limit of the periodic solutions $u_j$
to \eqref{NLS1} with $u_j|_{t = 0} = f_j$, where $f_j$ is as in \eqref{LP2}.
At this point, however, we do not know how to prove uniqueness in $C(\R; S^{1, 2}(\R))$.

\smallskip

\noi
(iii) It would be of interest to characterize limit periodic functions satisfying \eqref{LP3}.
 In Appendix \ref{SEC:4}, we present a brief discussion on  a sufficient condition for \eqref{LP3}.

\end{remark}

\begin{remark}\rm 

Let us compare Theorem  \ref{THM:1} and the result in \cite{BE}, when $k = 1$.
On the one hand,  the rate of approximation of $f$ by the periodization $f_j$ (i.e. the condition 
\eqref{LP3} in Theorem \ref{THM:1}) is more restrictive than that in \cite{BE}.
On the other hand, 
Theorem \ref{THM:1} holds with $s = 1$, 
while the result in \cite{BE} requires a higher regularity $s \geq 3$.

In Theorem \ref{THM:1}, 
we set $s = 1$.
Indeed, it is possible to lower the value of $s$
in view of global well-posedness of the defocusing periodic NLS
(in particular, for $k = 1$ and $2$; see \cite{BO, BO4}).
In this case, one needs to 
(i) control the growth of the $H^s$-norm of a solution
to the periodic problem
in the sprit of the $I$-method \cite{CKSTT}
and (ii) adjust the convergence rate \eqref{LP3}
to the regularity $s< 1$.
In our almost/limit periodic setting, 
it is important to control the $L^\infty$-norm (in $x$).
Hence, we need  $s > \frac{1}{2}$ in view of Sobolev embedding theorem.
\end{remark}

\begin{remark}\rm 
While we state and prove Theorem \ref{THM:1}
only for the defocusing case, 
the global existence result also holds  in the focusing case when $ k = 1$.
In this case, we need to replace \eqref{LP3} by 
\begin{align}
e^{L_{j+1}^{6 + \eps}} \| f _j - f\|_{S^{1, 2}(\R)} \leq B
\label{decay}
\end{align}

\noi
for all sufficiently large $j$.
See Remark \ref{REM:focusing}.

In the focusing case with $k \geq 2$, 
finite time blowup solutions are known 
to exist in $H^1(\R)$ and $H^1(\T)$.
In these settings, there are some criteria on such finite time blowup solutions
(such as negative energy \cite{OT2, OT1}).
Since a periodic function is in particular almost periodic,
 finite time blowup solutions in the periodic setting 
 \cite{OT2}
provide an instance of finite time blowup results
in the almost periodic setting (where initial data are periodic).
It would be interesting to provide 
a criterion for 
finite time blowup solutions to \eqref{NLS1}
in a generic (i.e.~non-periodic) almost periodic setting.

In \cite{Oh}, we studied 
the following NLS:
\begin{align}
i \dt u + \dx^2 u = \ld  |u|^{2k}
\label{NLS2}
\end{align}

\noi
 in a generic  almost periodic setting
and provided a criterion
for finite time blowup solutions,
depending only on the signs of 
the real and imaginary parts of $\ld$
and the mean value $M(f)$ of an almost periodic initial condition $f$.

\end{remark}

\begin{remark}\label{REM:Morrey} \rm

A version of Morrey's inequality states that
\begin{equation}
 |f(x) - f(x_0)| \leq C r^{1-\frac{d}{p}}\bigg(\int_{B(x_0, 2r)}|\nb f (y)|^p dy\bigg)^\frac{1}{p}
\label{Morrey1}
 \end{equation}

\noi
for all $x \in B(x_0, r)$ and  $d < p \leq \infty$.
Here, $B(x_0, r) \subset \R^d$ denotes the ball of radius $r$ centered at $x_0 \in \R^d$.
Then, given a limit periodic function $f \in S^{1, 2}(\R)$, 
it follows from \eqref{Morrey1}
 that $f \in C^{\frac 12}(\R)$.
Hence, the Fourier series of $f$ converges to $f$ uniformly.
See \cite[p.\,46]{B}.
In particular, by Lemma \ref{LEM:AP2}, we conclude that 
the function $f$ is indeed represented  by its Fourier series.

\end{remark}

\begin{remark}\label{REM:SP}\rm
The space $S^p$, $p\geq1$, of Stepanov's generalized almost periodic
functions
is precisely the closure of 
 the trigonometric polynomials
under the $S^p$-metric $d_{S^p}$ defined by 
\begin{align}
d_{S^p}(f, g):= \sup_{y \in \R} \bigg(\int_{y}^{y+1} |f(x) - g(x) |^p dx \bigg)^\frac{1}{p}.
\label{Stepanov1}
\end{align}

\noi
Note that this metric is induced by the $S^{s, p}$-norm with $s = 0$.
Unlike (uniformly) almost periodic functions, 
the $S^p$-generalized almost periodic functions
are  determined
up to sets of measure 0.
On the one hand, an almost periodic function
is uniformly continuous.
On the other hand, if $f \in S^p$ is uniformly continuous on $\R$,
then it is 
(uniformly) almost periodic.
See \cite[Theorem 6.16 on p.\,174]{C}.

\end{remark}

We conclude this introduction by stating a useful
lemma, 
allowing  us to extract a periodic component
from an almost periodic function. 
Let $f$ be a function on $\R$.
Given $n \in \mathbb{N}$
and $L > 0$, 
define the averaging operator $A_{n, L}$
by 
\begin{align}
 A_{n, L}[f](x) 
  : =  \frac{1}{n}\Big\{
f(x) + f(x + L) + \cdots + f\big(x + (n-1) L\big) 
\Big\}.
\label{ave0}
\end{align}

\noi
Then, we have the following convergence property
of the averaging operator
on almost periodic functions.

\begin{lemma}[Theorem on p.\,44 of \cite{B}]\label{LEM:ave}
Let $f$ be an almost periodic function with the Fourier series \eqref{AP3c}.
Then, for each $L >0$, 
the limit
\begin{align}
 f^{(L)}(x) 
 & := \lim_{n \to \infty} A_{n, L}[f](x) 
\label{ave}
\end{align}

\noi
exists uniformly in $x \in \R$.
Moreover, $f^{(L)}$
is a periodic function with period $L$
whose Fourier series consists of the terms 
of the Fourier series \eqref{AP3c} of $f$ which have period $L$.
Namely, 
\[ f^{(L)}(x) \sim \sum_{\o \in \Z / L} \ft f(\o) e^{2\pi i  \o x}.\]

\end{lemma}

\noi
In the following, 
we refer to $f^{(L)}$ as the periodization of $f$ (with period $L$).

\section{Sobolev spaces on a scaled torus and scaling invariance of NLS}
\label{SEC:2}

In this section,
 we briefly go over the basic definitions and properties 
of Sobolev spaces  
on a scaled torus $\T_\ld := \R /(\ld \Z)$,  
$\ld \geq 1$, 
along with the scaling symmetry of NLS \eqref{NLS1}.
Given a function $F$ on $\T_\ld$, we define 
its Fourier coefficient  by 
\begin{equation}
 \ft{F}\big(\tfrac n\ld\big) = \frac{1}{\ld} \int_{\T_\ld} F(x) e^{-2\pi i \frac{n}{\ld}x} dx, \qquad n \in \Z.
\label{Fourier0a}
 \end{equation}

\noi
We have the following Fourier inversion formula:
\begin{equation} F (x) =  
\sum_{ n \in \Z } \ft{F}\big(\tfrac{n}{\ld}\big) e^{2\pi i \frac{n}{\ld} x}
\label{Fourier0}
\end{equation}

\noi
and Plancherel's identity:
\[ \| F \|_{L^2(\T_\ld)} = \ld^\frac{1}{2}\| \ft F \|_{\l^2(\Z/\ld)}
= \ld^{\frac{1}{2}} \bigg( \sum_{n \in \Z } 
\big|\ft F\big(\tfrac{n}{\ld}\big)\big|^2\bigg)^\frac{1}{2}.
\]

\noi
 Note that the definition \eqref{Fourier0a}
 of the Fourier coefficient for periodic functions
 agrees  with the definition \eqref{AP2a} of the Fourier coefficient
for almost periodic functions.

Next, we define the homogeneous Sobolev spaces $\dot H^s(\T_\ld)$
and 
 the inhomogeneous Sobolev spaces $ H^s(\T_\ld)$
by the norms:
\begin{align}
\| F \|_{\dot H^s(\T_{\ld})} 
:= \ld^\frac{1}{2} \bigg( \sum_{n \in \Z } 
\big|\tfrac{n}{\ld}\big|^{2s} \big|\ft F\big(\tfrac{n}{\ld}\big)\big|^2\bigg)^\frac{1}{2}, 
\label{Sobolev0}\\
\| F \|_{ H^s(\T_{\ld})}
:= \ld^\frac{1}{2} \bigg( \sum_{n \in \Z } 
\big\langle\tfrac{n}{\ld}\big\rangle^{2s} \big|\ft F\big(\tfrac{n}{\ld}\big)\big|^2\bigg)^\frac{1}{2},
\label{Sobolev0a}
\end{align}

\noi
where $\jb{\, \cdot\, } = (1 + |\cdot|^2)^\frac{1}{2}$.

NLS \eqref{NLS1} on $\R$
enjoys several symmetries.
In particular, 
the scaling symmetry plays an essential role
in the proof of Theorem \ref{THM:1}.
The scaling symmetry states that 
if $u$ is a solution to  \eqref{NLS1}
with initial condition $\phi$, then
the scaled function $u^\ld$, defined by 
\begin{align}
u^{\ld} (t, x) = 
\mathfrak{S}_\ld[u](t, x) : = 
\ld^{-\frac{1}{k}} u(\ld^{-2} t, \ld^{-1} x),
\label{scaling1}
\end{align}
	
\noi
is also a solution to \eqref{NLS1}
with the scaled 
initial condition:
\begin{align}
f^{\ld} (x) = 
\mathfrak{S}_\ld[f](x) 
:= 
\ld^{-\frac{1}{k}} f( \ld^{-1} x).
\label{scaling2}
\end{align}

\noi
With a slight abuse of notation, 
we use $\mathfrak{S}_\ld$
to denote both
the dilation operator \eqref{scaling1} for functions in $t$ and $x$
and 
the dilation operator \eqref{scaling2}
for functions only  in $x$, depending on the context.

We conclude this section by discussing
 the effect of the scaling \eqref{scaling2}
 on different norms.
Let $f$ be a function on $\T$.
It follows from \eqref{Fourier0a}
and \eqref{scaling2} that 
the Fourier coefficient
of the scaled function $f^\ld := 
\mathfrak{S}_\ld[f] $ on $\T_{\ld}$ 
is given by 
\begin{align}
\ft{f^\ld}\big(\tfrac n\ld \big) = \ld^{-\frac{1}{k}}\ft{f}( n).
\label{scaling4}
\end{align}

\noi
Then, from \eqref{Sobolev0}, we have 
\begin{align}
\|f^\ld \|_{\dot H^s(\T_{\ld})} = \ld^{\frac 12 -s-\frac{1}{k}} 
\|f \|_{\dot H^s(\T)}.
\label{scaling3}
\end{align}

\noi
In particular, 
for  $\ld \geq 1$ and $ s \geq 0$, 
we have 
\begin{align}
 \|f \|_{ H^s(\T)} \les \ld^{- \frac 12 + s+\frac 1k}  \|f^\ld \|_{ H^s(\T_\ld)} 
\label{scaling5}
\end{align}

Let $F$ be  a function on $\T_\ld$.
Then, from \eqref{Fourier0} and \eqref{Sobolev0}, 
we have the following Sobolev embedding type estimate:
\begin{align}
\|F\|_{L^\infty(\T_\ld)} 
& 
\leq 
 |\ft{F}(0)|
+ C\ld^{s-\frac{1}{2}} \|F \|_{\dot H^s(\T_\ld)}, 
\label{Sobolev1}
\end{align}

\noi
as long as  $ s> \frac 12$.
Combining \eqref{scaling4},
\eqref{scaling3},   and \eqref{Sobolev1}, we obtain
\begin{align}
\|f^\ld\|_{L^\infty(\T_\ld)}
 \leq 
\ld^{-\frac{1}{k}} |\ft{f}(0)|
+ C
\ld^{-\frac{1}{k}}
 \|f \|_{\dot H^s(\T)}
\sim 
\ld^{-\frac{1}{k}}  \|f \|_{ H^s(\T)}
\label{scaling6}
\end{align}

\noi
for $f$ on $\T$ as long as  $ s> \frac 12$.

Lastly, we state a version of 
 Sobolev embedding on $\T_\ld$.
By Cauchy-Schwarz inequality along with a Riemann sum approximation, we have
\begin{align}
 \|F\|_{L^\infty(\T_\ld)} 
& \leq  \|\ft F\|_{\l^1(\Z/\ld)} 
\leq
 \bigg( \frac{1}{\ld} \sum_{n \in \Z } 
\frac{1}{(1+|\frac{n}{\ld}|^{2})^s} \bigg)^\frac{1}{2}
 \|F \|_{ H^s(\T_\ld)} \notag \\
& \les
 \|F \|_{ H^s(\T_\ld)}
\label{Sobolev2}
 \end{align}

\noi
for  $ s> \frac 12$.
Here, the implicit constants are independent of the period $\ld$.

\begin{remark} \rm
Note that  \eqref{Sobolev2} with 
\eqref{scaling3} only yields
\begin{equation*} 
\|f^\ld\|_{L^\infty(\T_\ld)}
 \les  \ld^{\frac{1}{2}- \frac 1k } \|f \|_{ H^s(\T)}
\end{equation*}

\noi
for $s > \frac 12$,
which 
is 
not as efficient as \eqref{scaling6}.
This is due to the fact that 
 the homogeneous Sobolev norms
act better than 
the inhomogeneous Sobolev norms with respect to the scaling.
In the following, we will use both \eqref{scaling6}
and \eqref{Sobolev2}.

\end{remark}

\section{
Global existence} 

In this section, we present the proof of Theorem \ref{THM:1}.
The proof is 
 based on an elementary scaling argument
 and the $H^1$-global well-posedness of the defocusing NLS in the periodic setting.
We first introduce some notations.
Given $j \in \mathbb{N}$,
we set 
\begin{equation}
\T_{j} : = \R / (L_j \Z),
\label{TJ}
\end{equation}

\noi
where $L_j =   j! /\o  $ 
as in \eqref{period}.
Given 
 a limit periodic function $f \in LP_\o(\R)$, 
let $f_j \in H^1(\T_j ) $ be 
the periodization of $f$ with period $L_j$ 
as in \eqref{LP2}. 
We assume that $f$ and $f_j$
satisfy the hypothesis in Theorem \ref{THM:1},
in particular, \eqref{LP3}.
In view of the $H^1$-global well-posedness 
of the defocusing NLS 
in the periodic setting, 
it follows that,
for each $j \in \mathbb{N}$, 
 there exists
a unique global solution 
 $u_j \in C(\R; H^1(\T_j))$ 
 to  \eqref{NLS1} with $u_j|_{t = 0} = f_j \in H^1(\T_j)$.

In the following, 
we show that $\{u_j\}_{j = 1}^\infty$
is a Cauchy sequence in $C(\R; L^\infty(\R))$
endowed with the compact-open topology (in $t$ with values in $L^\infty(\R)$).\footnote{
Recall that 
a sequence $\{w_j \}_{j = 1}^\infty \subset C(\R_t; L^\infty(\R_x))$ converges to $w$ in the compact-open topology 
if and only if, for every compact subset $K$ of $\R_t$, 
the sequence  $\{w_j (t) \}_{j = 1}^\infty$ converges  to $w(t) $ in $L^\infty(\R_x)$  uniformly in  $t \in K$.}
For this purpose, 
we perform two kinds of scalings
to $f_j$ and $u_j$.
For $j \in \mathbb{N}$, set
 \begin{align*} 
 g_{j} := \mathfrak{S}_{L_{j}}^{-1}[f_{j}] \qquad
 \text{and} 
\qquad  v_{j} :=  \mathfrak{S}_{L_{j}}^{-1}[u_{j}],
\end{align*}

\noi
where $\mathfrak{S}_\ld$ is the dilation operator defined in \eqref{scaling1} and \eqref{scaling2}.
We also set\footnote{Since $f_j $ and $u_j$ are periodic (in $x$) with period $L_{j}$, 
they are also periodic with period $L_{j+1} = (j+1) L_j$.} 
 \begin{align*} 
g^j :=  \mathfrak{S}_{L_{j+1}}^{-1}[f_j]
\qquad \text{and} \qquad  
v^j := \mathfrak{S}_{L_{j+1}}^{-1}[u_j].
\end{align*}

\noi
Note that $v_j$ and $v^j$ are global solutions to \eqref{NLS1} on 
the standard torus $\T$
with initial data $g_j$ and $g^j$, respectively.

We first establish an estimate on $v^j - v_{j+1}$.
By \eqref{scaling5} and \eqref{LP3}, we have 
\begin{align}
\| g^j \|_{H^1(\T)}
& \les  
L_{j+1}^{\frac{1}{2} +\frac{1}{k}}\| f_j \|_{H^1(\T_{j+1})}
\leq 
C(\|f\|_{S^{1, 2}(\R)})
L_{j+1}^{1+\frac 1k},  \label{F1}\\
\| g_{j+1} \|_{H^1(\T)}
& \les  
L_{j+1}^{\frac{1}{2} +\frac{1}{k}}\| f_{j+1} \|_{H^1(\T_{j+1})}
\leq 
C(\|f\|_{S^{1, 2}(\R)})
L_{j+1}^{1+\frac 1k},  \label{F2}
\end{align}

\noi
for sufficiently large $j \gg 1$.
We also have
\begin{align}
\| g^j - g_{j+1}\|_{H^1(\T)}
& \les L_{j+1}^{ \frac{1}{2} +\frac{1}{k}}\| f_j - f_{j+1}\|_{H^1(\T_{j+1})}\notag \\
& \leq  L_{j+1}^{ \frac{1}{2} +\frac{1}{k}}\big(\| f_j - f\|_{H^1(\T_{j+1})} + \|  f_{j+1} - f\|_{H^1(\T_{j+1})}\big) \notag \\
&  
 \les  
e^{- L_{j+1}^{2k + \eps}}
L_{j+1}^{1+\frac 1k}, 
\label{F3}
\end{align}

\noi
for sufficiently large $j \gg 1$.

Let $ 2 \leq p \leq \infty$.
Then, 
given a compact interval $I \subset \R$ with $|I| \geq 1$, 
we have the following  Gagliardo-Nirenberg inequality
\begin{align}
\| \phi \|_{L^p(I)}
 \les 
  \|  \phi \|_{L^2(I)}^{\frac{1}{2} + \frac{1}{p}}
 \|    \phi \|_{ H^1(I)}^{\frac{1}{2} - \frac{1}{p}}.
\label{Gag1}
\end{align}

\noi
See \cite{Brezis} for example.
By a simple scaling argument, 
we can choose
the implicit constant in \eqref{Gag1} to be  independent of $I$
with $|I| \geq 1$.
By \eqref{scaling2}, 
 \eqref{Gag1}, and \eqref{LP3}, we have 
\begin{align}
\|g^j \|_{L^{2k+2}(\T)}^{2k+2} 
& =  L_{j+1}^{1+\frac{2}{k}}\|f_j  \|_{L^{2k+2}(\T_{j+1})}^{2k+2} 
\leq  L_{j+1}^{2+\frac{2}{k}}
\sup_{\substack{I \subset \T_{j+1} \\ |I| = 1}}\|f_j  \|_{L^{2k+2}(I)}^{2k+2} \notag \\
& \les 
 L_{j+1}^{2+\frac{2}{k}}
\sup_{\substack{I \subset \T_{j+1} \\ |I| = 1}}
\|f_j \|_{L^2(I)}^{k+2} 
\|f_j \|_{ H^1(I)}^k 
\leq 
C(\|f\|_{S^{1, 2}(\R)})
 L_{j+1}^{2+\frac{2}{k}}
\label{F3a}
\end{align}

\noi
for sufficiently large $j \gg 1$.
Hence, by the conservation of the Hamiltonian
and the mass
with \eqref{F1} and \eqref{F3a},
we obtain  
\begin{align}
\| v^j(t) \|_{H^1(\T)} 
\leq 
C(\|f\|_{S^{1, 2}(\R)})
 L_{j+1}^{1+\frac{1}{k}}
\label{F3b}
\end{align}
	
\noi
for any $t \in \R$
and sufficiently large $j \gg 1$.
By a similar computation with \eqref{F2}, we also obtain
\begin{align}
\| v_{j+1}(t) \|_{H^1(\T)} 
\leq 
C(\|f\|_{S^{1, 2}(\R)})
 L_{j+1}^{1+\frac{1}{k}}
\label{F3c}
\end{align}

\noi
for any $t \in \R$ and sufficiently large $j \gg 1$.

Now, consider the Duhamel formulation of \eqref{NLS1} on $\T$:
\begin{align}
v^j(t) = \G_{g^j} v^j  (t) : = e^{i t \dx^2} g^j  - i \int_0^t e^{i (t-t') \dx^2}  |v^j|^{2k} v^j (t') dt'.
\label{NLS3}
\end{align}

\noi
By the unitarity of the linear propagator and the algebra property of $H^1(\T)$, 
one can easily show that the map $\G_{g^j}$ is a contraction on 
the ball of radius $2\|g^j\|_{H^1(\T)}$ in $C([-T_*, T_*]; H^1(\T))$
for some $T_* > 0$.
It follows from  \eqref{F1}
that   
this  standard fixed point argument
via the Duhamel formulation \eqref{NLS3}
yields the local time $T_*$ of existence,  
satisfying 
\begin{align*}
T_* \sim \|g^j\|_{H^1(\T)}^{-2k} \ges L_{j+1}^{-2(k+1)}
\end{align*}
	
\noi

\noi
for sufficiently large $j \gg 1$.
By repeating this argument with \eqref{F2}, 
we see that 
 the same argument holds 
even if we replace $g^j$ and $v^j$ with $g_{j+1}$ and $v_{j+1}$, 
respectively.
Note that, 
in view of the global-in-time control
\eqref{F3b} and \eqref{F3c}
of the $H^1$-norms of  the global solutions $v^j$ and $v_{j+1}$, 
we can iterate this local-in-time argument
indefinitely
for both $v^j$ and $v_{j+1}$
on time intervals of size 
\begin{align}
T_* 
 \sim L_{j+1}^{-2(k+1)}.
\label{F4}
\end{align}

\noi
Finally, consider the difference 
of the Duhamel formulations
for $v^j$ and $v_{j+1}$.
Then,  	
by iterating the local argument 
on intervals of size $T^*$ given by \eqref{F4}
and noting that 
the distance (in $H^1(\T)$) between $v^j$ and $v_{j+1}$
can grow at most by a fixed constant multiple on each 
of $O(\frac{T}{T^*})$ many such intervals, 
there exists $J_* \gg 1$ and $c > 0$ such that 
\begin{align}
\big\|  v^j - v_{j+1}\big\|_{C([-T, T]; H^1(\T))}
\les e^{  [c T L_{j+1}^{2(k+1)}]+1} \| g^j - g_{j+1}\|_{H^1(\T)}
\label{F5}
\end{align}

\noi
for all $T > 0$ and $j \geq J_*$.
Here, $[\tau]$ denotes the integer part of $\tau$.
Then, by undoing the scaling with \eqref{scaling6},  \eqref{F5}, and \eqref{F3},  we obtain 
\begin{align}
\| u_j -  u_{j+1}\|_{C([-T, T]; L^\infty(\R))}
&  = \| u_j - u_{j+1}\|_{C([-T, T]; L^\infty(\T_{j+1}))} \notag \\
& \les L_{j+1}^{-\frac1k} 
\|  v^j - v_{j+1}\|_{
C([-L_{j+1}^{-2}T, L_{j+1}^{-2}T]; H^1(\T))} \notag \\
&  \les e^{  [c T L_{j+1}^{2k}] + 1- L_{j+1}^{2k+\eps} } L_{j+1}
\too 0, 
\label{F6}
\end{align}

\noi
as $j \to \infty$, for each fixed $T> 0$.
Hence,  
from \eqref{F6} and \eqref{period}, 
we have 
\begin{align*}
\| u_j - u_{j'}\|_{C([-T, T];L^\infty(\R))}
& \leq \sum_{\l = j'}^{j - 1}e^{  - L_{\l+1}^{2k}} 
\sim e^{  - L_{j'+1}^{2k}}
\end{align*}

\noi
for $j \geq j' \gg 1$,
where the right-hand side converges to 0
as $j' \to \infty$.
Therefore, $\{u_j\}_{j = 1}^\infty$ converges in $C(\R; L^\infty(\R))$
with the compact-open topology, i.e.~for each $T> 0$, the convergence is uniform in $|t| \leq T$.
Denote the limit by $u\in C(\R; L^\infty(\R))$.
Then, we have $u|_{t = 0} = f$.
Moreover, 
the above convergence implies that  $u$ is a distributional solution to \eqref{NLS1}, i.e.~we have
\begin{align*}
\iint_{\R\times \R} u \Big( - i \dt \phi + \dx^2 \phi\Big) dx dt
=  
\iint_{\R\times \R} |u|^{2k} u  \cdot  \phi \, dx dt
\end{align*}

\noi
for any test function $\phi \in C^\infty_c(\R\times\R)$.

Since  $u_j (t) \in  H^1(\T_j)$ for each $t \in \R$, it follows that 
$u_j(t)$ is continuous (in $x$) for each $t \in \R$ and $j \in \mathbb{N}$.
Therefore, 
as a uniform limit of continuous periodic function $u_j(t)$, 
we conclude that $u(t)$ is 
limit periodic for each $t \in \R$.

We also claim that $u(t) \in LP_\o(\R)$ for each $t \in \R$.
Fix $t\in \R$.
Given $\eps > 0$, there exists $j \in \mathbb{N}$ such that 
\begin{equation}
\| u (t)- u_j (t) \|_{L^\infty(\R)} <  \frac{\eps}{3}.
\label{Y1}
\end{equation}

\noi
By  defining $u^{(L)}(t)$ and $u_j^{(L)}(t)$ be the periodizations of $u(t)$ and $u_j(t)$ as in  
\eqref{ave}, 
it follows from 
Lemma \ref{LEM:ave}
that there exists $N \in \mathbb{N}$
such that
\begin{align}
\sup_{x \in \R} 
\big|u_*^{(L)}(t, x)
- A_{n, L}[u_*](t,  x)
\big|
< \frac{\eps}{3}, 
\label{Y2}
\end{align}
	
\noi
for all $n \geq N$, where $u_* = u$ or $ u_j$.
Here, $A_{n, L}$ is the averaging operator defined  in \eqref{ave0}.
Then, from \eqref{Y1}
and \eqref{Y2} with \eqref{ave0}, 
we see that 
\begin{align}
\|u^{(L)}(t) - u^{(L)}_j(t)\|_{L^\infty(\R)} < \eps.
\label{Z1}
\end{align}

\noi
Now, let 
$L \in \R \setminus (\mathbb{Q}/\o)$.
Then, noting that $u_j(t)$ is periodic (in $x$) with period $L_j =  j!  /\o \in \mathbb Q/\o$, 
it follows from 
 Lemma \ref{LEM:ave} that  
 \begin{align}
 u^{(L)}_j (t) \equiv 0.
\label{Z2}
 \end{align}

\noi
Hence, from \eqref{Z1} and \eqref{Z2}, 
we have  
\[\|u^{(L)}(t)\|_{L^\infty(\R)} < \eps.\]

\noi
Since the choice of $\eps >0$ was arbitrary, we conclude that 
\[u^{(L)}(t)\equiv 0.\]

\noi
In particular, it follows from 
 Lemma \ref{LEM:ave} that $1/ L  \notin \s(u(t))$. 
Therefore, we have  
$\s(u(t)) \subset \o \cdot \mathbb{Q}$ for any $t \in \R$.
This proves that $u(t) \in LP_\o(\R)$
for each $t \in \R$.

Next,  we prove uniqueness  of the solution $u$ constructed above.
Given $\o \in \R \setminus\{0\}$, 
let $f \in LP_\o(\R)$ satisfy the hypothesis of Theorem \ref{THM:1}.
Then, it follows from Remark \ref{REM:Morrey}
that $f$ is given by 
its Fourier series \eqref{LP1}.
Set 
\[\oo := \{ r_m \o \}_{m = 1}^\infty \in \R^\mathbb{N},\]

\noi 
 where $r_m$ is as in \eqref{LP1}.
Then, 
uniqueness follows 
from the local well-posedness result  in 
$\A_{\oo}(\R)$
presented in  \cite{Oh}, 
once we show that 
 $u(t) \in  \A_{\oo}(\R)$
for each $t\in \R$.

\smallskip

\noi
$\bullet$
{\bf Case (a):} $t = 0$.

Let $A(j)$ be as in \eqref{LP2a}.
Then, we have 
\begin{equation}
\mathbb N = \bigcup_{j\in \mathbb N} A(j)
= A(1) \cup \bigcup_{j\in \mathbb N} \big\{A(j+1) \setminus A(j)\big\}
\label{U0}
\end{equation}

\noi
since $r_m \in \mathbb{Q}$ satisfies 
(i) $r_m  j! \in \Z$ for some $j = j(r_m) \in \mathbb{N}$
and (ii) $r_m  \wt j! \in \Z$ for all $\wt j \geq j  (r_m)$.

It follows from \eqref{Sobolev2} and \eqref{LP3}
that 
there exists $j_0 \in \mathbb {N}$
such that 
\begin{align}
\sum_{m \in A(j+1)\setminus A(j)}
|\ft f (r_m \o)|
& \les \|f_{j+1} - f_j \|_{H^1(\T_{j+1})}
\leq L_{j+1}^\frac{1}{2} \|f_{j+1} - f_j \|_{S^{1, 2}(\R)} \notag \\
& \leq L_{j+1}^\frac{1}{2} \big(\|f_{j+1} - f \|_{S^{1, 2}(\R)}
+ \|f - f_j \|_{S^{1, 2}(\R)}\big)\notag \\
& \les
e^{-  L^{2k}_{j+1}}
\label{U1}
\end{align}
	
\noi
for all $j \geq j_0$.
Similarly, we have
\begin{align}
\| f_{j_0}\|_{\A_{\oo}(\R)} 
& \les \| f_{j_0}\|_{H^1(\T_{j_0})} 
\leq L_{j_0}^{\frac{1}{2}}\big(\| f \|_{S^{1, 2}(\R)}
+ \| f_{j_0} - f \|_{S^{1, 2}(\R)}\big) \notag \\
& \les L_{j_0}^\frac{1}{2}\| f \|_{S^{1, 2}(\R)}
\label{U1a}
\end{align}

\noi
for sufficiently large $j_0 \gg 1$.
Then, from \eqref{U1} and \eqref{U1a} with \eqref{period}, 
we obtain
\begin{align}
\| f\|_{\A_{\oo}(\R)} 
&  = \sum_{m = 1}^\infty |\ft f (r_m \o)|
= \sum_{m \in A_{j_0}}|\ft f (r_m \o)| 
+ \sum_{j = j_0}^\infty
 \sum_{m \in A(j+1)\setminus A(j)}
|\ft f (r_m \o)|\notag\\
& \les \| f_{j_0}\|_{\mathcal{A}_{\oo}(\R)}
+ \sum_{j = j_0}^\infty e^{- L^{2k}_{j+1}}
\les L_{j_0}^\frac{1}{2}\|f\|_{S^{1, 2}(\R)}
+ e^{-  L^{2k}_{j_0+1}}
<\infty. \notag
\end{align}

\noi
Therefore, we conclude that $f \in A_{\oo}(\R)$.

\smallskip

\noi
$\bullet$
{\bf Case (b):} $t \ne 0$.

Fix $t \in \R \setminus \{0\}$.
Then, by slightly modifying the  computations in \eqref{F6}, we have 
\begin{align}
\| u_j  -  u_{j+1}  \|_{C([-T, T]; S^{1, 2}(\R))}
&  \leq  \| u_j  -  u_{j+1}  \|_{C([-T, T]; H^1(\T_{j+1}))} \notag \\
& \les L_{j+1}^{\frac 12 - \frac 1k } 
\|  v^j  -  v_{j+1} 
\|_{ C([-L^{-2}_{j+1}T, L^{-2}_{j+1}T];H^1(\T))}  \notag \\
&  \les e^{  [c T  L_{j+1}^{2k}] + 1- L_{j+1}^{2k+\eps} } L_{j+1}^{\frac 32 }
 \too 0, 
 \label{U2}
\end{align}

\noi
as $j \to \infty$,  for each $T> 0$.
In particular, $u_j(t)$ converges to $u(t)$ in $S^{1, 2}(\R)$,
uniformly
on the time interval $[-T, T]$ for each $T>0$.
Therefore, we conclude that $u\in C(\R; S^{1, 2}(\R))$.

Fix $j \in \mathbb{N}$.
Let $A_{n, L_j}$ be the averaging operator defined in \eqref{ave0}
and  $u^{(L_j)}(t)$ and $u_j^{(L_j)}(t)$ be the periodizations of $u(t)$ and $u_j(t)$ 
with period $L_j$   defined in \eqref{ave}.
By Lemma \ref{LEM:ave}, 
$A_{n, L_j} [u(t) - u_j(t)]$ converges uniformly (in $x$)
to $u^{(L_j)}(t) - u^{(L_j)}_j(t)$
 as $n \to \infty$.
In particular, 
given an interval $I \subset \R$ with $|I| = 1$, 
$A_{n, L_j} [u(t) - u_j(t)]$ converges $u^{(L_j)}(t) - u^{(L_j)}_j(t)$
in $L^2(I)$.
Moreover, 
we have 
\[\big\|A_{n, L_j} [u(t) - u_j(t)]\big\|_{H^1(I)}
\leq \|u(t) - u_j(t)\|_{S^{1, 2}(\R)}.\]

\noi
Namely, 
$\big\{A_{n, L_j} [u(t) - u_j(t)]\big\}_{n \in \mathbb{N}}$ is bounded in $H^1(I)$.
Therefore, 
$A_{n, L_j} [u(t) - u_j(t)]$ converges weakly in $H^1(I)$
as $n \to \infty$.\footnote{Suppose that $\{f_n\}_{n = 1}^\infty$
converges to $f$ in $L^2$ and is bounded in $H^1$.
Fix $\eps > 0$.
Given a test function $\phi \in H^{-1}$,
let $\phi_\eps  \in L^2$ such that 
$\| \phi - \phi_\eps\|_{H^{-1}} <  \eps/(2M)$, where $M = \sup_n \|f_n - f\|_{H^1}$.
Then, we have 
\begin{align*}
|\jb{ f_n - f, \phi}|
\leq 
\| f_n - f \|_{L^2}\| \phi_\eps\|_{L^2}
+ \tfrac{1}{2}\eps
< \eps
\end{align*}

\noi
for $n \geq N = N(\eps, \phi)$.
}
As a result, we obtain
\begin{align}
\| u^{(L_{j})}(t) -u_{j}^{(L_{j})}(t)\|_{S^{1, 2}(\R)} 
&  = 
\Big\| \lim_{n\to \infty} A_{n, L_j} [u(t) -u_{j}(t)]\Big\|_{S^{1, 2}(\R)}\notag \\
& \leq \liminf_{n \to \infty}
\big\|  A_{n, L_j} [u(t) -u_{j}(t)]\big\|_{S^{1, 2}(\R)} \notag \\
& \leq   \|u(t) - u_j(t)\|_{S^{1, 2}(\R)}. 
\label{U4}
 \end{align}

Since $u_j(t)$ is already periodic with period $L_j$, 
we have  $ u_j^{(L_j)}(t) = u_j(t)$.
Then, it follows from  the triangle inequality
with \eqref{U4} and  \eqref{U2}
that there exists $j_0\in \mathbb{N}$
such that 
\begin{align}
\| u^{(L_{j+1})} &  (t)  - u^{(L_{j})}(t) \|_{H^1(\T_{j+1})}  \notag \\
& \leq 
\big(\| u^{(L_{j+1})}(t) -u_{j+1}^{(L_{j+1})}(t)\|_{H^1(\T_{j+1})} \notag \\
& \hphantom{XXX} 
+ \| u_{j+1} (t) -  u_{j} (t) \|_{H^1(\T_{j+1})}
+ \| u_{j}^{(L_j)}(t) - u^{(L_{j})}(t)\|_{H^1(\T_{j+1})} \big) \notag\\
& \leq L_{j+1}^\frac{1}{2}
\big(\| u^{(L_{j+1})}(t) -u_{j+1}^{(L_{j+1})}(t)\|_{S^{1, 2}(\R)} \notag \\
& \hphantom{XXX} 
+ \| u_{j+1} (t) -  u_{j} (t) \|_{S^{1, 2}(\R))}
+ \| u_{j}^{(L_j)}(t) - u^{(L_{j})}(t)\|_{S^{1, 2}(\R)} \big) \notag\\
& \leq L_{j+1}^\frac{1}{2}
\big(\| u(t) -u_{j+1}(t)\|_{S^{1, 2}(\R)} \notag \\
& \hphantom{XXX} 
+ \| u_{j+1} (t) -  u_{j} (t) \|_{S^{1, 2}(\R))}
+ \| u_{j}(t) - u(t)\|_{S^{1, 2}(\R)} \big) \notag\\
& \les
e^{-   L^{2k}_{j+1}}
\label{U5}
\end{align}
	
\noi
for all $j \geq j_0$. 
Finally, proceeding as in \eqref{U1a}
with  \eqref{U5}, we obtain 
\begin{align*}
\| u(t) \|_{\A_{\oo}(\R)} 
&  = \sum_{m = 1}^\infty |\ft u (t, r_m \o)|  \notag \\
& = \sum_{m \in A_{j_0}}|\ft u (t, r_m \o)| 
+ \sum_{j = j_0}^\infty
 \sum_{m \in A(j+1)\setminus A(j)}
|\ft u (t, r_m \o)|\notag\\
& \les \| u^{(L_{j_0})}(t) \|_{H^1(\T_{j_0})}
+ \sum_{j = j_0}^\infty
\| u^{(L_{j+1})}(t) - u^{(L_{j})}(t) \|_{H^1(\T_{j+1})} \notag \\
& \leq L_{j_0}^\frac{1}{2}\| u(t) \|_{S^{1, 2}(\R)}
+ Ce^{- c L^{2k}_{j_0+1}}
<\infty.
\end{align*}

\noi
Therefore, we conclude that $u(t) \in A_{\oo}(\R)$
for each $t \in \R$.

Lastly, we present the proof of Theorem \ref{THM:1} (ii).
Given $J \in \mathbb{N}$ and $K > 0$, let 
 $\mathcal{B}^\o(J, K)$ be as in \eqref{BJ}.
Fix $T>0$ and $\eps > 0$.
Given $f, \wt f \in  \mathcal{B}^\o(J, K)$, 
let $u$ and $\wt u$ be the global solutions to \eqref{NLS1}
constructed above
with $f$ and $\wt f$ as initial data, respectively.
Denoting
  the periodizations   of $f$ and $\wt f$ with period $L_j$
  by  $f_j$ and  $\wt f_j$, 
we  denote by  $u_j$ and $\wt u_j$
the global solutions to \eqref{NLS1} on $\T_j$ with initial
data $f_j$ and $\wt f_j$, respectively.
Then, it follows from \eqref{U2}
that there exists $J_1 \in \mathbb{N}$
such that 
\begin{align}
\| u - u_j\|_{C([-T, T]; S^{1, 2}(\R))}
+ \| \wt u - \wt u_j\|_{C([-T, T]; S^{1, 2}(\R))}
< \frac{1}{2}\eps
\label{CT1}
\end{align}
	
\noi
for all $j \geq J_1$.

Let $v_j = \mathfrak{S}_{L_j}^{-1}[u_j]$
and $\wt v_j = \mathfrak{S}_{L_j}^{-1}[\wt u_j]$. 
By proceeding as in \eqref{U2}, we have 
\begin{align}
\| u_j - \wt u_j\|_{C([-T, T]; S^{1, 2}(\R))}
& \leq 
\| u_j - \wt u_j\|_{C([-T, T]; H^1(\T_j))}\notag \\
& \les L_j^{\frac 12 -\frac 1k} \| v_j - \wt v_j \|_{C([-L_j^{-2}T, L_j^{-2} T]; H^1(\T))}.
\label{CT2}
\end{align}
	
\noi	
Note that  $v_j$ and $\wt v_j$
satisfy the global $H^1$-bound \eqref{F3c}
(with $j+1$ replaced by $j$), where the constant depends only on $K$.
Hence, by iterating the local-in-time argument
over time intervals of size $T_* \sim L_j^{-2(k+1)}$, 
we obtain
\begin{align}
\| v_j - \wt v_j  & \|_{C([-L_j^{-2}T, L_j^{-2} T]; H^1(\T))}
 \leq e^{[cTL_j^{2k}] + 1}
\| v_j(0) - \wt v_j(0)\|_{H^1(\T)} \notag \\
& \leq e^{[cTL_j^{2k}] + 1} L_j^{\frac 12 + \frac 1k}
\| f_j  - \wt f_j \|_{H^1(\T_j)}
   \leq e^{[cTL_j^{2k}] + 1} L_j^{1 + \frac 1k}
\| f_j  - \wt f_j \|_{S^{1, 2}(\R)}\notag \\
&  \leq e^{[cTL_j^{2k}] + 1} L_j^{1 + \frac 1k}
\big( \| f_j  - f \|_{S^{1, 2}(\R)} 
\notag\\
& 
\hphantom{XXXXX}
+ \| f  - \wt f \|_{S^{1, 2}(\R)}
+ \| \wt f  - \wt f_j \|_{S^{1, 2}(\R)}\big).
\label{CT3}
\end{align}

\noi
It follows from \eqref{CT2}
and \eqref{CT3} with \eqref{LP3}
that there exists $J_2 \in \mathbb{N}$
such that 
\begin{align}
\| u_j - \wt u_j\|_{C([-T, T]; S^{1, 2}(\R))}
&  \leq e^{[cTL_j^{2k}] + 1} L_j^{\frac 32}
 \| f  - \wt f \|_{S^{1, 2}(\R)}  + \frac{1}{4} \eps
\label{CT4}
\end{align}

\noi
for all $j \geq J_2$.

Finally, 
letting $j_* = \max(J_1, J_2)$, 
it follows from \eqref{CT1}
and \eqref{CT4} that 
\begin{align*}
\| u - \wt u \|_{C([-T, T]; S^{1, 2}(\R))}
&  \leq e^{[cTL_{j_*}^{2k}] + 1} L_{j_*}^{\frac 32}
 \| f  - \wt f \|_{S^{1, 2}(\R)}  + \frac{3}{4} \eps
 < \eps, 
\end{align*}
	
\noi	
where the last inequality holds
as long as we have  
\[  \| f  - \wt f \|_{S^{1, 2}(\R)} < \dl  = \dl (\eps, j_*) \ll 1.\]

\noi
Note that by  choosing  $J_1, J_2 \geq J$, 
we can make sure that 
 they  do not depend on 
a particular choice of functions in $\mathcal{B}^\o(J, K)$.
This proves continuous dependence of the flow on $\mathcal{B}^\o(J, K)$.

\begin{remark}\label{REM:focusing} \rm 
In the following, we briefly discuss how Theorem \ref{THM:1}
also follows in the focusing case if $k = 1$.
Let $\phi \in H^1(\T)$.
Then,  it follows from  \eqref{Gag1}
that there exists $C_0> 0$ such that 
\begin{align*}
\| \phi\|_{H^1(\T)}^2 
\les M[\phi] + H[\phi]
+ C_0 \big(M[\phi]\big)^3
\end{align*}

\noi
In particular, by the conservation of the Hamiltonian and the mass
along with \eqref{scaling3}, 
we have 	
\begin{align}
\| v^j(t) \|_{H^1(\T)}
&  \leq \Big(M[g^j] + H[g^j]
+ C_0 \big(M[g^j]\big)^3\Big)^\frac{1}{2} \notag \\
& \leq 
C(\|f\|_{S^{1, 2}(\R)})
 L_{j+1}^{3 }
\end{align}

\noi
in place of \eqref{F3b}.	
A simliar computation 	
shows that \eqref{F3c} also holds with 	
$ L_{j+1}^{3}$.
As a result, we obtain $T^* \sim 
 L_{j+1}^{-6}$ instead of \eqref{F4}.
Then, it is easy to see that the rest of the proof of Theorem \ref{THM:1}
goes 
through with small modifications,
as long as \eqref{decay} holds.

\end{remark}

\appendix

\section{On the condition \eqref{LP3}}
\label{SEC:4}

In this appendix, we briefly investigate a meaning of the condition \eqref{LP3}
on the rate of convergence of $f_j$ to a limit periodic initial condition $f$.
Given a quasi-periodic function $f$ of the form \eqref{AP4a}, one can talk about a decay of
the Fourier coefficients 
of the form:
$|\ft f (\pmb{\o}\cdot \pmb{n})| \les | \pmb{n}|^{-\g}$
and $|\ft f (\pmb{\o}\cdot \pmb{n})| \les \exp( - \kappa | \pmb{n}|^{\theta})$
for $ \pmb{n} \in \Z^N$.
For a generic almost periodic function, however, 
such conditions do not make sense since $| \pmb{n}|  =\infty $
for $ \pmb{n} \in \Z^\mathbb{N}$, 
unless $\pmb{n} = (n_1, n_2, \dots)$ has a finite support.
Since our limit periodic initial condition $f$ in Theorem \ref{THM:1} is not quasi-periodic, 
it does not seem appropriate or at least seems non-trivial to characterize the condition \eqref{LP3}
in terms of a decay  of the Fourier coefficients only in $|\pmb{n}|$.
We instead consider a sufficient condition for \eqref{LP3}
and discuss a decay of the Fourier coefficients in $n_j$  (see \eqref{APP1} below
for the definition of $n_j$) and in $L_j$ (and hence in $j$)
in the following.

Let $A(j)$ be as in \eqref{LP2a} with the understanding that $A(0) = \emptyset$.
With \eqref{LP1} and Remark \ref{REM:Morrey}, we have
\begin{align}
 f(x) & = \sum_{m = 1}^\infty \ft f(r_m \o) e^{ 2\pi i r_m \o x}
= \sum_{j = 1}^\infty \sum_{m \in A(j) \setminus A(j-1)} \ft f(r_m \o) e^{ 2\pi i r_m \o x}. \notag
\end{align}

\noi
Here, each summand in the $j$-summation is given by $f_j - f_{j-1}$ defined in \eqref{LP2}, 
and thus is periodic
with period $L_j = j!/\o$.
With  $n_j = r_m j! $, we have
\begin{align}
 f(x) & = \sum_{j = 1}^\infty \sum_{\substack{n_j \in \Z \\ j\, \nmid  \,   n_j}}
\ft f\big(\tfrac{n_j}{L_j}\big) e^{ 2\pi i \frac{n_j}{L_j} x}.
\label{APP1}
 \end{align}

\noi
Note that we have 
\[\ft f\big(\tfrac{n_j}{L_j}\big) = \ft f_j\big(\tfrac{n_j}{L_j}\big),\]

\noi
where the Fourier transform on the right-hand side is 
that for periodic functions with period $L_j$ discussed
in  Section \ref{SEC:2}.
With  $\pmb{n} = (n_1, n_2, \dots)$
and $\pmb{\o} = (L_1^{-1}, L_2^{-1}, \dots)$, 
we see that \eqref{APP1} is 
an analogous formulation to \eqref{AP4a}  in our limit periodic setting,
showing that we have infinite many summations over $\Z$ unlike the quasi-periodic setting.

%

In the remaining part, we consider the following sufficient
condition for \eqref{LP3}:
\begin{align}
 \|  f_{j} -f_{j-1}  \|_{S^{1, 2}(\R)} \leq C e^{- L_{j}^{2k + \eps}}
\label{APP2}
\end{align}

\noi
for all sufficiently large $j$.
This in turn is guaranteed if we have 
\begin{align}
 \|  f_{j} -f_{j-1}  \|_{H^1(\T_j)} \leq C e^{- L_{j}^{2k + \eps}}, 
\label{APP3}
\end{align}

\noi
where $\T_j$ as in \eqref{TJ}.
Hence, by letting  $F_j = f_j - f_{j-1}$ denote the difference
of consecutive periodizations $f_{j}$ and $f_{j-1}$ (with the understanding that $f_0 \equiv 0$), 
we see that the condition \eqref{LP3} is satisfied
if (i) $f = \sum_{j = 1}^\infty F_j$
and 
(ii) the $H^1(\T_j)$-norms 
of the $L_j$-periodic functions $F_j$
decay at 
a rate $C e^{- L_{j}^{2k + \eps}}$
for all sufficiently large $j$.

Lastly, note that the condition \eqref{APP3} is essentially 
necessary for \eqref{APP2}.
Indeed, 
 it follows from \eqref{APP1} and \eqref{APP2}
that 
\begin{align*}
C e^{- L_{j}^{2k + \eps}}
& \geq 
 \|  f_{j} -f_{j-1}  \|_{S^{1, 2}(\R)}
\geq L^{-\frac 12}_j \|  f_{j} -f_{j-1}  \|_{H^1(\T_j)}.
\end{align*}

\noi
Hence, from \eqref{Sobolev0a}, we must have 
\begin{align}
\bigg(\sum_{\substack{n_j \in \Z \\ j\, \nmid  \,   n_j}}
\big\langle \tfrac{n_j}{L_j} \big\rangle^2 
\big|\ft f\big(\tfrac{n_j}{L_j}\big) \big|^2\bigg)^\frac{1}{2}
\leq 
C' e^{- L_{j}^{2k + \frac{\eps}{2}}}
\label{last}
\end{align}

\noi
for all sufficiently large $j$, if \eqref{APP2} holds.
The condition \eqref{last} states
that the Fourier coefficients
$\ft f\big(\tfrac{n_j}{L_j}\big)$
must decay  polynomially in $n_j$ (in an average sense)
besides the very fast decay in $j$
for each $L_j$-periodic component.

\begin{acknowledgment} \rm
The author would like to thank
the Hausdorff Research Institute for Mathematics
for its generous hospitality
during the author's stay in 
the  trimester program: ``Harmonic Analysis and Partial Differential Equations'',
where a part of this manuscript was prepared.
The author  is also grateful to the anonymous referees
for their helpful comments.
\end{acknowledgment}


\begin{thebibliography}{99}		

\bibitem{B}

A.~Besicovitch, {\it Almost periodic functions,} Dover Publications, Inc., New York, 1955. xiii+180 pp.

\bibitem{Bohr}
H.~Bohr,  {\it  Zur theorie der fast periodischen funktionen. I. Eine verallgemeinerung der theorie der fourierreihen,}  Acta Math. 45 (1925), no. 1, 29--127.

\bibitem{BO}
J.~Bourgain, {\it Fourier transform restriction phenomena for certain lattice subsets and applications to nonlinear evolution equations. I. Schr\"odinger equations,} Geom. Funct. Anal. 3 (1993), no. 2, 107--156.



\bibitem{BO4}
J.~Bourgain, {\it A remark on normal forms and the ``$I$-method'' for periodic NLS,}  J. Anal. Math.  94  (2004), 125--157.



\bibitem{BE}
A.~Boutet de Monvel, I.~Egorova,
{\it On solutions of nonlinear Schr\"odinger equations with Cantor-type spectrum,}
J. Anal. Math. 72 (1997), 1--20. 


\bibitem{Brezis}
H.~Brezis, 
{\it Functional analysis}, Sobolev spaces and partial differential equations. Universitext. Springer, New York, 2011. xiv+599 pp.


\bibitem{CAZ}
T.~Cazenave, 
{\it Semilinear Schr\"odinger equations,}
Courant Lecture Notes in Mathematics, 10. New York University, Courant Institute of Mathematical Sciences, New York; American Mathematical Society, Providence, RI, 2003. xiv+323 pp. 


\bibitem{CKSTT}
J.~Colliander, M.~Keel, G.~Staffilani, H.~Takaoka, T.~Tao, 
{\it Sharp Global Well-Posedness for KdV and Modified KdV
on $\mathbb{R}$ and $\mathbb{T}$,}
J. Amer. Math. Soc. 16 (2003), no. 3, 705--749.




\bibitem{C}
C.~Corduneanu, {\it Almost periodic functions,} 
With the collaboration of N. Gheorghiu and V. Barbu. Translated from the Romanian by Gitta Bernstein and Eugene Tomer. Interscience Tracts in Pure and Applied Mathematics, No. 22. Interscience Publishers [John Wiley \& Sons], New York-London-Sydney, 1968. x+237 pp.

\bibitem{DG}
D.~Damanik, M.~Goldstein, 
{\it On the existence and uniqueness of global solutions for the KdV equation with quasi-periodic initial data}, 
arXiv:1212.2674v2 [math.AP]. 

\bibitem{E}
I.~Egorova, {\it The Cauchy problem for the KdV equation with almost periodic initial data whose spectrum is nowhere dense,} Spectral operator theory and related topics, 181--208, Adv. Soviet Math., 19, Amer. Math. Soc., Providence, RI, 1994.

\bibitem{GV}
J.~Ginibre, G.~Velo,
{\it On a class of nonlinear Schr\"odinger equations. III. Special theories in dimensions 1, 2 and 3,}
Ann. Inst. H. Poincar\'e Sect. A (N.S.) 28 (1978), no. 3, 287--316. 



\bibitem{K}
Y.~Katznelson, {\it An introduction to harmonic analysis,} Third edition. Cambridge Mathematical Library. Cambridge University Press, Cambridge, 2004. xviii+314 pp. 





\bibitem{OT2}
T.~Ogawa, Y.~Tsutsumi, 
{\it Blow-up of solutions for the nonlinear Schr\"odinger equation with quartic potential and periodic boundary condition,} Functional-analytic methods for partial differential equations (Tokyo, 1989), 236--251, Lecture Notes in Math., 1450, Springer, Berlin, 1990. 


\bibitem{OT1}
T.~Ogawa, Y.~Tsutsumi, 
{\it Blow-up of $H^1$ solutions for the one-dimensional nonlinear Schr\"odinger equation with critical power nonlinearity,} Proc. Amer. Math. Soc. 111 (1991), no. 2, 487--496.


\bibitem{Oh}
T.~Oh, {\it On nonlinear Schr\"odinger equations with almost periodic initial data,}
 to appear in SIAM J. Math. Anal. 

\bibitem{TAO}
T.~Tao, {\it Nonlinear dispersive equations. Local and global analysis,} CBMS Regional Conference Series in Mathematics, 106. Published for the Conference Board of the Mathematical Sciences, Washington, DC; by the American Mathematical Society, Providence, RI, 2006. xvi+373 pp. 

\bibitem{Tsugawa}
K.~Tsugawa, {\it Local well-posedness of the KdV equation with quasi-periodic initial data,} SIAM J. Math. Anal. 44 (2012), no. 5, 3412--3428.


\bibitem{Tsutsumi}
Y.~Tsutsumi, {\it $L^2$-solutions for nonlinear Schr\"odinger equations and nonlinear groups,} Funkcial. Ekvac. 30 (1987), no. 1, 115--125.



\end{thebibliography}
\end{document}